\documentclass[11pt]{amsart}

\usepackage{amsmath,amssymb,amscd,mathrsfs}

\theoremstyle{definition}
\newtheorem{thm}{Theorem}[section]
\newtheorem{lem}[thm]{Lemma}
\newtheorem{prp}[thm]{Proposition}
\newtheorem{dfn}[thm]{Definition}
\newtheorem{cor}[thm]{Corollary}

\newtheorem{cnj}[thm]{Conjecture}

\newtheorem{rmk}[thm]{Remark}

\newtheorem{exa}[thm]{Example}

\newenvironment{pff}{{\em Proof:}}{\QED}

\newcommand{\beq}{\begin{equation}}
\newcommand{\eeq}{\end{equation}}
\newcommand{\beqr}{\begin{eqnarray*}}
\newcommand{\eeqr}{\end{eqnarray*}}

\newcommand{\bit}{\begin{itemize}}
\newcommand{\eit}{\end{itemize}}

\newcommand{\af}{\alpha}
\newcommand{\bt}{\beta}

\newcommand{\dt}{\delta}

\newcommand{\ld}{\lambda}

\newcommand{\ph}{\varphi}

\newcommand{\vd}{\vartheta}

\newcommand{\Ld}{\Lambda}
\newcommand{\Sm}{\Sigma}

\newcommand{\Om}{\Omega}

\newcommand{\Cg}{\mathcal C}
\newcommand{\Dg}{\mathcal D}

\newcommand{\uld}{{\underline{\ld}}}
\newcommand{\umu}{{\underline{\mu}}}

\newcommand{\ui}{{\underline{i}}}
\newcommand{\uj}{{\underline{j}}}

\newcommand{\uLd}{{\underline{\Ld}}}

\newcommand{\hg}{{\hat{g}}}

\newcommand{\Cart}{{\mathrm{Cart}}}

\newcommand{\Z}{{\mathbb{Z}}}

\newcommand{\C}{{\mathbb{C}}}
\newcommand{\N}{{\mathbb{N}}}

\renewcommand{\hg}{{\hat{\mathfrak{g}}}}
\newcommand{\Sl}{{\mathfrak{sl}}}

\newcommand{\g}{{\mathfrak{g}}}

\newcommand{\U}{{U_{\Z}}}
\newcommand{\V}{{V_{\Z}}}
\renewcommand{\H}{{\mathcal{H}}}

\newcommand{\diag}{{\mathrm{diag}}}

\newcommand{\Par}{{\mathrm{Par}}}

\newcommand{\la}{{\langle}}
\newcommand{\ra}{{\rangle}}

\newcommand{\andeqn}{\,\,\,\,\,\, {\mbox{and}} \,\,\,\,\,\,}
\newcommand{\QED}{$\Box$}

\begin{document}
\footnotetext{Version of September 24, 2008}
\begin{center}\textbf{Cartan Invariants of Symmetric Groups and Iwahori-Hecke Algebras}\end{center}

\bigskip

\begin{center}
Christine Bessenrodt

Institut f\"ur Algebra, Zahlentheorie und Diskrete Mathematik

Leibniz Universit\"at Hannover, Welfengarten 1

D-30167 Hannover, Germany

\verb"bessen@math.uni-hannover.de"

\end{center}

\bigskip

\begin{center}

David Hill

University of California, Berkeley

Berkeley, California

U.S.A.

\verb"dhill1@math.berkeley.edu"
\end{center}

\begin{abstract}
K\"{u}lshammer, Olsson and Robinson conjectured that a certain set of numbers determined the invariant factors of the $\ell$-Cartan matrix for $S_n$ (equivalently, the invariant factors of the Cartan matrix for the Iwahori-Hecke algebra $\mathcal{H}_n(q)$, where $q$ is a primitive $\ell$th root of unity). We call these invariant factors Cartan invariants. 

In a previous paper, the second author calculated these Cartan invariants when $\ell=p^r$, $p$ prime, and $r\leq p$ and went on to conjecture that the formulae should hold for all $r$. Another result was obtained, which is surprising and
counterintuitive from a block theoretic point of view. Namely, given the prime decomposition $\ell=p_1^{r_1}\cdots p_k^{r_k}$, the Cartan matrix of an $\ell$-block of $S_n$ is a product of Cartan matrices associated to $p_i^{r_i}$-blocks of $S_n$. In particular, the invariant factors of the Cartan matrix associated to an $\ell$-block of $S_n$ can be recovered from the Cartan matrices associated to the $p_i^{r_i}$-blocks.

In this paper, we formulate an explicit combinatorial determination of the Cartan invariants of $S_n$--not only for the full Cartan matrix, \emph{but for an individual block}. We collect evidence for this conjecture, by showing that the formulae predict the correct determinant of the $\ell$-Cartan matrix. We then go on to show that Hill's conjecture implies the conjecture of KOR.
\end{abstract}

\section{Introduction}

The theory of generalized blocks of symmetric groups was initiated by
K\"{u}lshammer, Olsson and Robinson in \cite{kor}.
Using character-theoretic methods,  they showed that
many invariants of the usual block theory of symmetric groups over a field of characteristic $p$
do not depend on $p$ being a prime. This led the authors to define `$\ell$-blocks' of
symmetric groups and a related $\ell$-modular representation theory.
They defined an appropriate analogue of the Cartan matrix associated to $S_n$ for this theory and even
conjectured that a certain set of numbers determined the invariant factors of this matrix
\cite[Conjecture 6.4]{kor}. In a related paper \cite{bo}, Bessenrodt and Olsson
conjectured a formula for the determinant of the Cartan matrix.

Using a new method developed in \cite{llt, ar, gr, kl}, Brundan and Kleschev \cite{bk}
calculated an explicit
formula for the determinant of the Cartan matrix of a block of the Iwahori-Hecke algebra,
$\H_n$, with
parameter $q$ a primitive $\ell$th root of unity.
Donkin \cite{d} showed that there is a direct link between
$\ell$-blocks of $S_n$ and blocks of $\H_n$.
In particular, their respective Cartan matrices have the same
determinant and invariant factors. Using this,
together with the results of \cite{bk} and \cite{bo},
K\"{u}lshammer, Olsson and Robinson \cite{kor}
verified the formula conjectured by Bessenrodt and Olsson \cite{bo}
(see also the remarks at the end of \cite{bo}).
It should also be noted that in \cite{bos},
Bessenrodt, Olsson and Stanley obtained a more elementary proof of the formula
for the determinant of the full Cartan matrix.

In \cite{h}, Hill investigated the invariant factors of the Cartan matrix
associated to an individual block of
$\H_n$ using the methods developed in \cite{bk}.
When $\ell=p^r$ is a power of a prime satisfying $r\leq p$
these numbers were computed (see \cite[Theorem 1.3]{h}).
Moreover, he conjectured that the same
formula held for arbitrary $r$.

In \cite{h} also another result was obtained, which is surprising and
counterintuitive from a block theoretic point of view. Namely, given the prime
decomposition $\ell=p_1^{r_1}\cdots p_k^{r_k}$, the Cartan matrix of an
$\ell$-block of $S_n$ is a product of Cartan matrices associated to
$p_i^{r_i}$-blocks of $S_n$. In particular, the invariant factors of the Cartan
matrix associated to an $\ell$-block of $S_n$ can be recovered from the Cartan
matrices associated to the $p_i^{r_i}$-blocks (see \cite[Theorem 1.1, 1.2]{h}).

We want to emphasize that -- going beyond the conjecture in \cite{kor} -- we
conjecture here an explicit combinatorial determination of the invariants not
only for the full $\ell$-Cartan matrix of~$S_n$, but even for the Cartan
matrices of the $\ell$-blocks, see Conjecture~\ref{conj:block Cartan}.
In our context this is a very natural refinement.
In principle, it should also be possible to obtain a block version of the
conjecture in \cite{kor} by using \cite[Theorem 6.1]{kor} and methods
similar to the ones applied in \cite{bo}; an explicit combinatorial
sorting of the invariants given by \cite{kor} into blocks has not been described so far,
though.

In the remainder of the paper, we collect evidence for this conjecture (and the
conjecture in \cite{h}, respectively), by showing that the formulae predict the
correct determinant for the $\ell$-Cartan matrices. We then go on to show that
the conjecture in \cite{h} implies the conjecture in \cite{kor}. Thus, in
particular, the results in \cite{h} mentioned above imply that the latter
conjecture holds when any prime divisor $p$ in $\ell$ occurs in $\ell$ with
exponent at most~$p$. For the convenience of the reader, we also calculate
\cite[Examples 6.5, 6.6, 6.7]{kor} using our methods (see Examples \ref{example
1}-\ref{example 3}).

We would like to point out that a lot of the machinery required for the proofs
in this article has already been developed in \cite{bo}. We find this striking,
and hope that the exposition here will help to elucidate the relationship
between this new approach to the representation theory of symmetric groups and
the classical block and character theoretic methods.

\bigskip

\section{Background and preliminaries}

\subsection{Kac-Moody Algebras, Iwahori-Hecke algebras and Cartan Matrices}

In this section we give a brief description of the connection between highest
weight representations of Kac-Moody algebras and the representation ring of
Iwahori-Hecke algebras. We refer the reader to \cite{bk} for details in this case,
and \cite{kl} for the general theory.

Let $\hg$ be the affine Kac-Moody algebra of type $A_{\ell-1}^{(1)}$,
working always over the field $\C$ of complex numbers. Let
$e_0,\ldots,e_{\ell-1}$, and $f_0,\ldots,f_{\ell-1}$ be the Chevalley
generators of $\hg$, and $\tau:U(\hg)\rightarrow U(\hg)$ the
Chevalley anti-involution defined by $\tau(e_i)=f_i$
($i=0,\ldots,\ell-1$). We are interested in the basic representation
$V=V(\Ld_0)$ of $\hg$. It is the irreducible highest weight
representation with highest weight satisfying $\Ld_0(h_i)=\dt_{i0}$
($h_i=[e_i,f_i]$), see \cite{kc}. Fix a nonzero highest weight vector
$v_+\in V$. The Shapovalov form $(\cdot,\cdot)_S:V\times
V\rightarrow\C$ is the unique Hermitian form on $V$ satisfying
$(v_+,v_+)_S=1$ and $(xv,v')_S=(v,\tau(x)v')_S$ for all $x\in U(\g)$
and $v,v'\in V$. The weights of $V$ are
of the form $w\Ld_0-d\dt$, where $w$ is an element of the affine Weyl
group, $d\in\Z_{\geq0}$, and $\dt$ is the null root (see \cite[ $\S12.6$]{kc}).
Let $\U$ be the Kostant-Tits $\Z$-subalgebra of
$U(\hg)$ generated by the divided powers $e_i^{(n)}:=e_i^n/n!$ and $f_i^{(n)}:=f_i^n/n!$
($0\leq i\leq \ell-1$, $n\geq 1$) in the Chevalley generators. Define
$\V=\U v_+$. The Shapovalov form restricts to a symmetric bilinear
form $(\cdot,\cdot)_S:\V\times\V\rightarrow\Z$.

The lattice $\V$ is related to the Iwahori-Hecke algebra,
$\H_n=\H_n(q)$, of the symmetric group $S_n$ over an algebraically
closed field $F$ ($\mathrm{Char}F=p\geq0$) and with finite quantum
characteristic $\ell$. The quantum characteristic of $\H_n$ is
defined to be the number
\[
\ell=\min\{\,e\geq 2\,|\,1+q+\cdots+q^{e-1}=0\,\}
\]
if it exists, and $\ell=\infty$ otherwise.
For finite $\ell$, $q=1$ implies $\ell=p$ and $\H_n=FS_n$ and $q\neq1$ implies $q$ is a primitive
$\ell$th root of unity. The algebra $\H_n$ is not semisimple. The
simple $\H_n$-modules are labeled by the set $\Par_{\ell}^*(n)$ of
$\ell$-regular partitions (see section~2.2),
and the same is true for their projective
covers (i.e., the projective indecomposable modules). The main problem
is to describe the composition multiplicities $[P_\ld:L_\mu]$ of
the simple module $L_\mu$ inside the projective cover $P_\ld$ of
$L_\ld$, $\ld,\mu\in\Par_{\ell}^*(n)$.

Let $K_n=K(\H_n)$ be the Grothendieck group of the category of
finitely generated projective $\H_n$-modules. The Cartan pairing
$(\cdot,\cdot)_C:K_n\times K_n\rightarrow\Z$ is defined on the
projective indecomposable modules by $(P_\ld,P_\mu)_C=[P_\mu:L_\ld]$.
The Grothendieck group $K_n$ decomposes into blocks, and two
irreducibles are in the same block if, and only if, the partitions
labeling them have the same $\ell$-core and $\ell$-weight, see
\cite{ma}, and the blocks of $K_n$ are orthogonal with respect to the
Cartan pairing. The \emph{Cartan matrix}
\[
C_{\ell}(n):=([P_\mu:L_\ld])_{\ld,\mu\in\Par_{\ell}^*(n)}
\]
is the Gram matrix of this form. The matrix $C_{\ell}(n)$ is block
diagonal with blocks corresponding to the blocks of $K_n$.

Now, by \cite{ar,gr,kl}, we have
\[
\V\cong \bigoplus_{n\geq 0}K_n=:K
\]
as $\U$-modules, with the action of the $e_i^{(n)}$ (resp. $f_i^{(n)}$) are described in terms of certain
restriction functors (resp. induction functors). Under this isomorphism, the Shapovalov form corresponds
to the Cartan
pairing and the $(w\Ld_0-d\dt)$-weight space of $\V$ corresponds to
the block of $K$ with $\ell$-core associated to $w\Ld_0$ and
$\ell$-weight $d$, see \cite[$\S5.3$]{llt}  for details.

\subsection{Partitions and Multipartitions} Let $\Par(d)$ be the set
of all partitions of $d$, and $p(d)=|\Par(d)|$.
Given an integer $\ell\geq 1$, we say that
$\ld\in\Par(d)$ is $\ell$-class regular if no part of $\ld$ is
divisible by $\ell$, and we say $\ld$ is $\ell$-regular if no part of
$\ld$ is repeated $\ell$ (or more) times. Let $\Par_{\ell}(d)$ (resp.
$\Par_{\ell}^*(d)$) denote the set of all $\ell$-class regular (resp.
$\ell$-regular) partitions of~$d$, and $p_{\ell}(d)=|\Par_{\ell}(d)|$
($p_{\ell}^*(d)=|\Par_{\ell}^*(d)|$, resp.).
Finally, define
\[
\Par=\bigcup_{d\geq0}\Par(d),\;\;\;\Par_{\ell}=
\bigcup_{d\geq0}\Par_{\ell}(d),\andeqn\Par_{\ell}^*=\bigcup_{d\geq0}\Par_{\ell}^*(d).
\]
It is well known that the generating function
$P(q)=\sum_{d\geq0}p(d)q^d$ is given by
\[
P(q)=\prod_{i\geq 1}\frac{1}{1-q^i}.
\]
The generating function $P_{\ell}(q)$ for the numbers $p_{\ell}(d)$ is
\begin{eqnarray}\label{formula:P_l}
P_{\ell}(q)=\frac{P(q)}{P(q^\ell)}.
\end{eqnarray}
There is a bijection $G:\Par_{\ell}(d)\rightarrow\Par_{\ell}^*(d)$ known
as the Glashier bijection \cite{gl}. Hence, $P_{\ell}(q)$ is the
generating function for $p_{\ell}^*(d)$ as well.

Define the set of $\ell$-multipartitions of $d$ to be
\[
M_{\ell}(d)=\{\,\uld=(\ld^{(1)},\ldots,\ld^{(\ell)})\,|\,\ld^{(i)}\in\Par(d_i)\mbox{
for }1\leq i\leq\ell\mbox{ and }d_1+\cdots+d_{\ell}=d\}.
\]
and set $M_{\ell}=\bigcup_{d\geq0}M_{\ell}(d)$. The generating function
for the numbers $k(\ell,d)=|M_{\ell}(d)|$, i.e.,
$\sum_{d\geq 0} k(\ell, d)q^d$, is just $P(q)^\ell$.

\subsection{Divisors, the Total Length Function, and Cartan Matrices}
In this section, we review some facts about generating functions that
will be used in calculations below; these are mostly contained in \cite{bo}.

First, observe that the
generating function for the number of divisors of an integer $d$ is
\[
T(q)=\sum_{i\geq1}\frac{q^i}{1-q^i}.
\]
For a partition $\ld$, let $l(\ld)$ denote its length, i.e., the number of
its (non-zero) parts.
Then $l(d)=\sum_{\ld\in\Par(d)}l(\ld)$ is the total length function,
with corresponding generating function
$L(q):=\sum_{d\geq0}l(d)q^d$.
This is related to the number of divisors of $d$ by the equation
(see \cite[Proposition 2.1]{bo})
\begin{eqnarray}\label{formula:LPT}
L(q)=P(q)T(q).
\end{eqnarray}

More generally, these functions have $\ell$-class regular versions. Indeed,
\begin{eqnarray}\label{formula:T}
T(q)=T(q^\ell)+T_{\ell}(q)
\end{eqnarray}
where $T_{\ell}(q)$ is the generating function for the number of
divisors of $d$ which are not divisible by~$\ell$. Let
$
l_{\ell}(d)=\sum_{\ld\in \Par_{\ell}(d)}l(\ld)
$
be the total length function for the class $\ell$-regular partitions,
and $L_{\ell}(q):=\sum_{d\geq0}l_{\ell}(d)q^d$.
Then one has (see \cite[Proposition 2.2]{bo})
\begin{eqnarray}\label{formula:l-LPT}
L_{\ell}(q)=P_{\ell}(q)T_{\ell}(q).
\end{eqnarray}
Now, one easily concludes (see \cite[Corollary 2.3]{bo})
\begin{eqnarray}\label{formula:L-dec}
L(q)=P_{\ell}(q)L(q^\ell)+P(q^\ell)L_{\ell}(q).
\end{eqnarray}

We now turn our attention to some facts about the determinant of the
Cartan matrix $C_{\ell}(n)$ for the Iwahori-Hecke algebra $\H_n$ with
quantum characteristic $\ell$. As explained in the introduction, this
matrix encodes the composition multiplicities of simple modules
inside projective indecomposable modules.
Let $\ell^{c_{\ell}(n)}=\det C_{\ell}(n)$
and $\Cg_{\ell}(q):=\sum_{n\geq0}c_{\ell}(n)q^n$.
>From the work by Brundan and Kleshchev \cite{bk} and Bessenrodt and Olsson \cite[Theorem 3.3]{bo}
one obtains the following result  (see the remarks at the end of \cite{bo}):

\begin{thm}\label{thm:Cartan-det}
We have
$$\Cg_{\ell}(q)=P_{\ell}(q)T(q^\ell)\:.$$
\end{thm}

For later purposes we note that from this result
and (\ref{formula:P_l}) we may immediately deduce
the following, perhaps surprising, reduction formula:
\begin{cor}\label{cor:Cartan-det}
Let $a,b\in \N$. Then
$$\Cg_{ab}(q)=P_a(q)\Cg_b(q^a)\:.$$
\end{cor}

Since we usually work with a block version of this determinant, we
explain how to reconstruct the full determinant from this data.
Blocks of $C_{\ell}(n)$ are labeled by $\ell$-cores in $\Par(n-\ell w)$,
$0\leq w\leq\lfloor\frac{n}{\ell}\rfloor$. Here, $w$ is the
associated $\ell$-weight of the block. Let $d_{\ell}^0(n)$ be the number of
$\ell$-cores in $\Par(n)$. Then the corresponding generating function
$\Dg_{\ell}^0(q):=\sum_{n\geq0}d_{\ell}^0(n)q^n$
is given by
\begin{eqnarray}\label{formula:coregenfct}
\Dg_{\ell}^0(q)=\frac{P(q)}{P(q^\ell)^\ell}
\end{eqnarray}
(see \cite{o}). Let $\ell^{b_{\ell}(w)}$ be the determinant of a Cartan matrix
of a block of weight $w$, and let $B_{\ell}(q)$ be the associated
generating function for the numbers $b_\ell(w)$.
The connection between the generating functions for the full Cartan determinant
and the block Cartan determinant is then given by
\begin{eqnarray}\label{full-and-block}
\Cg_{\ell}(q)=B_{\ell}(q)\Dg_{\ell}^0(q).
\end{eqnarray}
Now \cite[ Theorem 3.4]{bo}, together with the remarks at the end of \cite{bo}
and (\ref{formula:LPT}) give
\begin{eqnarray}\label{block-det}
B_{\ell}(q)=P(q)^{\ell-1}T(q)=P(q)^{\ell -2} L(q).
\end{eqnarray}

\section{The Invariants}

The invariant factors of the Cartan matrix for $\H_n$ are determined
by the Shapovalov form on the lattice $\V$ as described in the
introduction. To explain the structure of $\V$, consider the simple
finite dimensional Lie algebra $\g$ (over $\C$), with (Lie) Cartan matrix
$A=(a_{ij})_{i,j=1}^{\ell-1}$, simple roots
$\{\af_1,\ldots,\af_{\ell-1}\}$, and root system
$Q=\bigoplus_{i}\Z\af_i$. Then, as $\U$-modules,
\[
V_\Z\cong\Z[Q]\otimes\uLd
\]
where $\Z[Q]$ is the group algebra of $Q$ and where
\[
\uLd=\bigotimes_{i=1}^{\ell-1}\Ld^{(i)},\andeqn
\Ld^{(i)}=\lim_{\substack{\leftarrow\\k}}\Z[x_1^{(i)},\ldots,x_k^{(i)}]^{S_k}
\]
is the ring of symmetric functions in the variables \emph{colored} by
$i$, see \cite[Theorem 4.5]{bk} and \cite[$\S3$]{h} for details.

Identifying $\V$ via this isomorphism, the highest weight vector is
then $e^0\otimes 1\in\V$, and a basis for $\V$ is given by
\[
\{\,e^\af\otimes h_\uld\,|\,\af\in Q,\,\uld\in M_{\ell-1}(d)\},
\]
where $h_\uld=\prod_{i=1}^{\ell-1}h_{\ld^{(i)}}(x^{(i)})$ if $\uld=(\ld^{(1)},\ldots,\ld^{(\ell-1)})$.

It was shown in \cite[Lemma 4.1]{bk}, that in this basis the Shapovalov form is given by
\[
(e^\af\otimes h_\uld,e^\bt\otimes h_\umu)_S=\dt_{\af\bt}\la h_\uld,h_\umu\ra_S,
\]
where $\la\cdot,\cdot\ra_S$ is the Shapovalov form on $\uLd$.
To describe this form more explicitly, let
\[
X_{A,d}=(\la m_\uld,h_\umu\ra_S)_{\uld,\umu\in M_{\ell-1}(d)}\andeqn X_A=\bigoplus_{d\geq0}X_{A,d}.
\]
The matrix $X_{A,d}$ can be regarded as a linear
transformation via the mapping
\[
\ph_A:\uLd_d\rightarrow\uLd_d,\;\;\;\ph_A(h_\uld)=\sum_{\umu\in
M_{\ell-1}(d)}(X_A)_{\uld\umu}h_\umu
\]
for all $\uld\in M_{\ell-1}(d)$. Let
\[
\Cart(A,d)=\uLd_d \big/\ph_A\left(\uLd_d\right)
\]
denote the corresponding finite group.
Note that $X_{A,d}$ has the same invariant factors as a submatrix of
$C_{\ell}(n)$ corresponding to an $\ell$-block of weight~$d$,
and that these also give the orders
of the cyclic factors of the finite abelian group $\Cart(A,d)$.

We now describe the matrix $X_A$. Let $\ld\in\Par$ and let
\[
\Om(\ld)=\{\ui=(i_1,\ldots,i_{l(\ld)})|1\leq i_k\leq\ell-1\mbox{ for
}1\leq k\leq l(\ld)\mbox{ and }i_j\leq i_{j+1}\mbox{ if
}\ld_j=\ld_{j+1}\}.
\]
Given a multipartition $\uld=(\ld^{(1)},\ldots,\ld^{(\ell-1)})\in
M_{\ell-1}(d)$, associate a pair $(\ld,\ui)$, where $\ld\in\Par(d)$
and $\ui\in\Om(\ld)$ as follows. Write $\uld$ as a single partition
$\ld=(\ld_1\geq\ld_2\geq\cdots)$. Each part $\ld_k$ belongs to some
$\ld^{(i_k)}$. Hence, we obtain a sequence
$\ui=(i_1,i_2,\ldots)\in\Om(\ld)$ by following the rule that $i_j\leq
i_{j+1}$ if $\ld_j=\ld_{j+1}$. The map $\uld\mapsto(\ld,\ui)$ is a
bijection, see \cite[Notation 3.1 and 3.2]{h} for details.
For any integer $r$, let $m_r(\ld)$ denote the multiplicity of $r$ as a part
of~$\ld$, and set
\begin{eqnarray*}
z_{\ld}=\prod_{r\geq 1}r^{m_r(\ld)}\cdot m_r(\ld)!\:.
\end{eqnarray*}
In this
notation, the Shapovalov form on $\uLd$ is defined on the power sum
symmetric functions by
\[
\la p_\uld,p_\umu\ra_S=\dt_{\ld\mu}a_{i_1j_1}a_{i_2j_2}\cdots a_{i_{\ell-1}j_{\ell-1}}z_\ld
\]
where $\uld\mapsto(\ld,\ui)$ and $\umu\mapsto(\ld,\uj)$
(see \cite[VI $\S10$]{m}).

Given a $k\times k$ matrix $Y=(y_{ij})_{i,j=1}^k$, define its $m$th
symmetric power $S^m(Y)$ to be the matrix with rows (resp. columns)
labeled by $m$-tuples $(i_1\leq i_2\leq\cdots\leq i_m)=:\ui$, and
with $(\ui,\uj)$ entry equal to
\[
y_{i_1j_1}y_{i_2j_2}\cdots y_{i_mj_m}.
\]
Define
\[
B_{A,d}=\bigoplus_{\ld\in\Par(d)}S^{m_1(\ld)}(A)\otimes\cdots\otimes
S^{m_d(\ld)}(A),\;\;\;B_A=\bigoplus_{d\geq 0}B_{A,d}.
\]
Notice that the rows and columns of the matrix $B_{A,d}$ are naturally
labeled by $M_{\ell-1}(d)$ via the bijection
$\uld\leftrightarrow(\ld,\ui)$ above.
Let $\Ld$ denote the (uncolored) ring of symmetric functions, and
$\Ld_d$ the its $d$th graded component (i.e., the span of all
symmetric functions of homogeneous degree $d$).
Let $M(p,m)$ be the transformation matrix between the power sum
basis and the monomial basis of $\Ld$, i.e., given by
$p_\ld=\sum_{\mu\in\Par}M(p,m)_{\ld\mu}m_\mu$.
We have, by \cite[$\S3$]{h},
\begin{eqnarray}\label{B_A}
X_A=(M(p,m)^{\otimes\ell-1})^{-1}B_A M(p,m)^{\otimes\ell-1}
\end{eqnarray}
where $M(p,m)^{\otimes\ell-1}$ is the matrix given by
$p_\uld=\sum_{\umu\in M_{\ell-1}} M(p,m)_{\uld\umu}m_\umu$.

Define a bilinear
pairing $\la\cdot,\cdot\ra_{\ell}:\Ld\times\Ld\rightarrow\Z$ on the
power sum symmetric functions by
\[
\la p_\ld,p_\mu\ra_{\ell}=\dt_{\ld\mu}\ell^{l(\ld)}z_\ld\:.
\]
Define the matrix $X_{\ell,d}=(\la
m_\ld,h_\mu\ra_{\ell})_{\ld,\mu\in\Par(d)}$, and
$X_{\ell}=\bigoplus_{d\geq0}X_{\ell,d}$. The matrix $X_{\ell}$ resembles
the matrix $X_A$ (they even agree when $\ell=2$). In fact, the matrix
$X_{\ell}$ is constructed in the same manner as $X_A$, with the matrix
$A$ replaced by the $1\times1$ matrix $(\ell)$. Indeed, let
\begin{eqnarray}\label{B}
B_{\ell,d}=\diag\{\ell^{l(\ld)}|\ld\in\Par(d)\}\andeqn
B=\bigoplus_{d\geq0}B_{\ell,d}.
\end{eqnarray}
We have (see \cite[$\S3$]{h})
\begin{eqnarray}\label{X}
X_{\ell}=M(p,m)^{-1}B_{\ell} M(p,m)\:.
\end{eqnarray}
As above, the matrix
$X_{\ell,d}$ can be regarded as a linear transformation via the map
$\ph_{\ell}:\Ld\rightarrow\Ld$,
$\ph_{\ell}(h_\ld)=\sum_\mu(X_{\ell,d})_{\ld\mu}h_\mu$. Let
$\Cart(\ell,d)=\Ld_d/\ph_{\ell}(\Ld_d)$ denote the corresponding finite
group.

Finally, we relate the matrix $X_{\ell}$ to $X_A$. To this end, recall
that the Smith normal form $\Sm(X)$ of a matrix $X$ is a diagonal
matrix with entries equal to the elementary divisors of $X$. Let $U$
and $V$ be unimodular matrices (i.e., integer matrices of determinant
$\pm1$) transforming the Lie Cartan matrix $A$ to its Smith form, i.e.,
\[
UAV=\Sm(A)=\diag\{1,1,\ldots,1,\ell\}.
\]
Define the matrix $B_U$ by the formulae
\[
B_{U,d}=\bigoplus_{\ld\in\Par(d)}S^{m_1(\ld)}(U)\otimes\cdots\otimes
S^{m_d(\ld)}(U),\\;\;\;B_U=\bigoplus_{d\geq 0}B_{U,d}
\]
and similarly for $B_V$.
Then, by \cite[Proposition 3.3]{h},
the matrices $X_U=M(p,m)^{-1}B_UM(p,m)$ and $X_V=M(p,m)^{-1}B_VM(p,m)$
are unimodular. Moreover, the matrix
\begin{eqnarray}\label{reduction of X}
\nonumber X_{U,d}X_{A,d}X_{V,d}&=&(I\otimes X_{\ell})_d\\
        &=&\bigoplus_{0\leq s\leq d}I_{d-s}\otimes X_{\ell,d}
\end{eqnarray}
where $I$ is the identity matrix on $\Ld^{\otimes\ell-2}$ (resp.
$I_{d-s}$ is the identity on the degree $d-s$ component of
$\Ld^{\otimes\ell-2}$). It follows that the matrix $I_{d-s}$ has rows
and columns labeled by $M_{\ell-2}(d-s)$. Finally, we observe that
this matrix has the same invariant factors as $X_{A,d}$. Hence, we
have the following:

\begin{thm}\label{inv factor multiplicities}(\cite[Theorem 1.1]{h})
Let $b_{1,s},\ldots,b_{h,s}$ be the invariant factors of
$X_{\ell,s}$ ($h=p(s)$), so $\Cart(\ell,s)$ is a finite abelian
group with cyclic factors of these orders. The finite abelian group
$\Cart(A,d)$ is a direct sum of $k(\ell-2,d-s)$
copies of $\Cart(\ell,s)$ for each $0\leq s\leq d$.
\end{thm}

In particular, the
coefficient of $q^{d-s}$ in the generating series $P(q)^{\ell-2}$
gives the number of cyclic factors
$\Z/b_{i,s}\Z$ that each $b_{i,s}$ contributes
to $\Cart(A,d)$.

\begin{rmk}\label{error} In \cite{h}, the statement of the theorem
above was (slightly) incorrect. The theorem stated that the invariant
factors of $X_{A,d}$ were the diagonal entries of the matrix
\[
\bigoplus_{0\leq s\leq d}I_{d-s}\otimes \Sm(X_{\ell,s}).
\]
In general, these are not equal to the diagonal entries of
\[
\Sm(X_{A,d})=\Sm\left(\bigoplus_{0\leq s\leq d}I_{d-s}\otimes
\Sm(X_{\ell,s})\right)
\]
unless $\ell$ is a power of a prime (though there is an easy
algorithm for going from the first set of diagonal entries to the
second in any specific case).
\end{rmk}

\begin{thm}\label{splitting}(\cite[Theorem 1.2]{h}) Let $a,b\in\Z_{\geq 2}$.
Then, $X_{ab,d}=X_{a,d}X_{b,d}$.
Moreover, if $(a,b)=1$, then
$\Sm(X_{ab,d})=\Sm(X_{a,d})\Sm(X_{b,d})$.
\end{thm}

\begin{pff} Using equations (\ref{B}) and (\ref{X}),
it is easy to see that $X_{ab,d}=X_{a,s}X_{b,d}$.
The second statement follows immediately by \cite[Theorem II.15]{n} since
\[
(\det X_{a,d},\det X_{b,d})=1.
\]
\end{pff}

Therefore, it is enough to compute the invariant factors of the
matrices $X_{p^r,d}=(X_{p,d})^r$ for every prime $p$ and $r\geq 1$.
In \cite{h}, these numbers were calculated in the case when $r\leq
p$. Indeed, let $\nu_p:\Z\rightarrow\Z_{\geq0}$ be the $p$-adic
valuation map (i.e., $\nu_p(n)=k$ if $n=p^kq$ with $(q,p)=1$). For a
positive integer $a$, define the number
\begin{eqnarray}\label{formula: defect}
d_p(a)=\sum_{j\geq 1}\left\lfloor\frac{a}{p^j}\right\rfloor\:.
\end{eqnarray}
Note that by the Legendre formula, it is known that
$p^{d_p(a)}=(a!)_p$ is the $p$-part of $a!$.
Let $\ld=(1^{m_1(\ld)}2^{m_2(\ld)}\ldots)$ and define
the $p$-defect of $\ld$ to be
\begin{eqnarray}\label{formula: ld-defect}
d_p(\ld)=\sum_{n\geq 1}d_p(m_n(\ld))
\end{eqnarray}
By a well-known result of Brauer on the invariants of the $p$-Cartan
matrix of a finite group, together with \cite{d}, we have
\begin{eqnarray}\label{eqn: Cartan-defect}
c_p(n)=\sum_{\ld\in\Par_p(n)}d_p(\ld)\:.
\end{eqnarray}

Furthermore, we define for $\ld=(1^{m_1(\ld)}2^{m_2(\ld)}\ldots)$:
\begin{eqnarray}\label{eqn: invariant factors}
\vd_{p^r}(\ld)=\prod_{\substack{n\geq1\\0\leq\nu_p(n)<r}}
p^{(r-\nu_p(n))m_{n}(\ld)+d_p(m_{n}(\ld))}.
\end{eqnarray}
Defining, for any integers $\ell,k\geq 1$, $\ell_k=\ell/(\ell,k)$,
we may then also write the numbers $\vd_{p^r}(\ld)$ in the form
\begin{eqnarray}\label{eqn: invariant factors 2}
\vd_{p^r}(\ld)=\prod_{\substack{n\geq1\\0\leq\nu_p(n)<r}} {(p^r)_n}^{m_n(\ld)}(m_n(\ld)!)_p.
\end{eqnarray}
Note that this clearly implies that
\[
\vd_{p^r}(1^d)=p^{rd}(d!)_p
\]
is the unique largest $p$-power among all $\vd_{p^r}(\ld)$, $\ld \in \Par(d)$.
\medskip

The following was proved in \cite[Theorem 1.3]{h}:
\begin{thm}\label{thm:main}
Let $r\leq p$. Then, the invariant factors of $X_{p^r,d}$ are the numbers
\[
\vd_{p^r}(\ld)\;,\ld\in\Par(d).
\]
In particular, $\vd_{p^r}(1^d)=p^{rd}d!_p$ is the unique largest invariant factor
of $X_{p^r,d}$.
\end{thm}

In the composite case the following definition is crucial.
\begin{dfn}\label{graded invariant factors}
Let $\ell=\prod_{i=1}^rp_i^{r_i}$ be the prime decomposition of $\ell$.
For $\ld\in\Par(d)$ we set
\[
\vd_{\ell}(\ld)=\prod_{i=1}^r\vd_{p_i^{r_i}}(\ld).
\]
Then, the \emph{graded invariant factors} of $X_{\ell,d}$ are defined to be
the numbers
\[
\,\vd_{\ell}(\ld)\, ,\ld\in\Par(d)\,.
\]
If $\ld\in\Par(d)$, we say that $\vd_{\ell}(\ld)$ has \emph{degree}~$d$.
\end{dfn}

\begin{cnj}\label{conj}
For any $\ell$, the finite group $\Cart(\ell,d)$
is a direct product of cyclic groups with orders given by the graded invariant
factors of $X_{\ell,d}$.
\end{cnj}

Observe that this conjecture is true under the assumption that
$\ell=\prod_{i=1}^rp_i^{r_i}$ satisfies $r_i\leq p_i$ for all $i$. In
the remaining sections we build evidence for this conjecture.
But first, we determine  the multiplicity of a graded invariant
factor in $C_{\ell}(n)$, and compute an example.

\begin{thm}\label{total multiplicity}
Let $m_{\ell}^n(d)$ be the coefficient of $q^{n-\ell d}$ in the
generating series
\[
\frac{P_{\ell}(q)}{P(q^\ell)}.
\]
Then, each $\ld \in \Par(d)$, $d\leq \left\lfloor \frac n \ell \right\rfloor$,
contributes $m_{\ell}^n(d)$
graded invariant factors $\vd_{\ell}(\ld)$ to $C_{\ell}(n)$.
\end{thm}

\begin{pff} The number of blocks at $n$ of weight $w$ is the
number of $\ell$-cores in $\Par(n-\ell w)$. This is
the coefficient of $q^{n-\ell w}$ in the generating series
$P(q)/P(q^\ell)^\ell$ by equation~(\ref{formula:coregenfct}).
Therefore, by Theorem~\ref{inv factor multiplicities}
and equation~(\ref{formula:P_l}),
each $\ld \in \Par(d)$
contributes the number $\vd_{\ell}(\ld)$ as a graded invariant factor of $C_{\ell}(n)$ $m_{\ld}$
times, where $m_\ld$ is the coefficient of $q^{n-\ell d}$ in
the generating series
\[
\frac{P(q)}{P(q^\ell)^\ell}P(q^\ell)^{\ell-2}=\frac{P_{\ell}(q)}{P(q^\ell)}\,.
\]
\end{pff}

\begin{exa}\label{example 1} Let $n=8$, $\ell=4$. We calculate the (graded) invariant factors of the
(principal) $\ell$-block
of weight 2 as given by the $\vd_4(\ld)$, where $\ld$ runs over all partitions of $d\leq2$. These are
$\vd_4(\emptyset)=1$,$\vd_4((1))=4$, $\vd_4((2))=2$, and $\vd_4((1^2))=32$ of degrees 0, 1, 2 and 2,
respectively. Finally, their respective multiplicities are $m_4^8(0)=5$, $m_4^8(1)=2$, $m_4^8(2)=1$, and
$m_4^8(2)=1$.

In summary, the graded invariant factors are $32^1,4^2,2^1,1^5$; here the exponents denote multiplicity.
Observe that  these graded invariants  coincide
with the numbers in \cite[Example 6.5]{kor}.
\end{exa}

\begin{exa}\label{example 2}Let $n=18$, $\ell=6$.
We give the graded invariant factors of $C_6(18)$ in the table below,
writing the number of contributions in each degree
towards the invariants of $C_6(18)$ in the form
\[
\sum_w(\mbox{mult. of a graded inv. factor in a block of weight
}w)\times(\mbox{\# of $\ell$-cores of }n-\ell w).
\]

\begin{table}[h]
\begin{tabular}{|c|l|c|}\hline Degree $d$&Graded Invariant Factors&Multiplicity in degree $d$\\
\hline0&1&$40\times1+14\times5+4\times20+1\times32=222$\\
1&6&$14\times1+4\times5+1\times20=54$\\
2&3, 72&$4\times1+1\times5=9$\\
3&2, 18, 1296&$1\times1=1$\\\hline
\end{tabular}

\vspace{.25in}

\caption{Graded Invariant Factors of $C_6(18)$}
\end{table}

Therefore, the graded invariant factors of $C_6(18)$ are:
$1^{222}$, $2^2$, $3^{9}$, $6^{54}$, $18^1$, $72^9$, $1296^1$;
the exponents denote multiplicity.
\end{exa}

\begin{exa}\label{example 3} Let $n=24$, $\ell=6$.
We give the graded invariant factors of $C_6(24)$ in the table below,
writing the number of contributions in each degree
towards the invariants of $C_6(24)$ in the form
\[
\sum_w(\mbox{mult. of a graded inv. factor in a block of weight
}w)\times(\mbox{\# of $\ell$-cores of }n-\ell w).
\]


\begin{table}[h]
\begin{tabular}{|c|l|c|}\hline Degree $d$&Graded Invariant Factors&Multiplicity in degree $d$\\
\hline0&1&$105\times 1+40\times 5+14\times 20+4\times32+1\times38=751$\\
1&6&$40\times1+14\times 5+4\times20+1\times32=222$\\
2&3, 72&$14\times1+4\times5+1\times20=54$\\
3&2, 18, 1296&$4\times1+1\times5=9$\\
4&3, 9, 12, 216, 31104&$1\times1=1$\\ \hline
\end{tabular}

\vspace{.25in}

\caption{Graded Invariant Factors of $C_6(24)$}
\end{table}

Therefore, the graded invariant factors of $C_6(24)$ are:
$1^{751}$, $2^9$, $3^{54+1}$, $6^{222}$, $9$, $12$, $18^9$,
$72^{54}$, $216$, $1296^9$, $31104$; the
exponents denote multiplicity.\\
The graded invariants of the Cartan matrix of the (principal) $6$-block of weight~$4$
are given by
$1^{105}$, $2^4$, $3^{14+1}$, $6^{40}$, $9$, $12$, $18^4$,
$72^{14}$, $216$, $1296^4$, $31104$.

It is instructive
to observe that most of the calculations required for this example are already
done in Example~\ref{example 2}.
\end{exa}

Again, note that  the graded invariants computed above for the full Cartan matrix
are precisely the numbers  $r_\ell(\ld)$, $\ld\in \Par_{\ell}(n)$,
in \cite[Examples 6.6, 6.7]{kor}. This is true in general, as
will be proved in section~\ref{kor}.

\section{The Determinant}

In this section we show that
$\prod_{\ld\in\Par(d)}\vd_{p^r}(\ld)=\det X_{p^r,d}$.
Because of equations~(\ref{B}) and (\ref{X}), we know that
\[
\det X_{p^r,d}=(p^{r})^{l(d)}\,.
\]
Therefore, it is enough to prove the following:

\begin{prp}
\[
\sum_{\ld\in\Par(d)}\log_p\vd_{p^r}(\ld)=r\,l(d).
\]
\end{prp}

\begin{pff}
Let $\ld=(1^{m_1(\ld)}2^{m_2(\ld)}\ldots)$. Recall from (\ref{eqn: invariant factors}) that
\begin{eqnarray*}
\log_p\vd_{p^r}(\ld)&=&\sum_{\substack{n\geq
1\\0\leq\nu_p(n)<r}}(r-\nu_p(\ld))m_n(\ld)+d_p(m_n(\ld))\\
    &=&\sum_{\substack{n\geq
    1\\(n,p)=1}}\sum_{i=0}^{r-1}(r-i)m_{p^in}(\ld)+d_p(m_{p^in}(\ld)).
\end{eqnarray*}
Now, given $\ld\in\Par(d)$, let $\ld=\ld^{(0)}+p\ld^{(1)}+\cdots
p^{N}\ld^{(N)}$ be the $p$-adic decomposition of $\ld$ (i.e.,
$\ld^{(i)}$ is a $p$-class regular partition, for all $i$). Then,
\begin{eqnarray*}
\log_p\vd_{p^r}(\ld)&=&
\sum_{n\geq 1}\sum_{i=0}^{r-1}(r-i)m_{n}(\ld^{(i)})+d_p(m_{n}(\ld^{(i)}))\\
    &=&\sum_{i=0}^{r-1}(r-i)l(\ld^{(i)})+d_p(\ld^{(i)})\:.
\end{eqnarray*}
Hence,
\begin{eqnarray*}
\sum_{\ld\in\Par(d)}\log_p\vd_{p^r}(\ld)&=&\sum_{\ld\in\Par(d)}\sum_{i=0}^{r-1}(r-i)l(\ld^{(i)})
+d_p(\ld^{(i)})\\
    &=&\sum_{\ld\in\Par(d)}\left(\sum_{i=0}^{r-2}(r-1-i)l(\ld^{(i)})+d_p(\ld^{(i)})\right)+
    \sum_{\ld\in\Par(d)}\left(\sum_{i=0}^{r-1}l(\ld^{(i)})+d_p(\ld^{(r-1)})\right).
\end{eqnarray*}
Therefore,
\begin{eqnarray}\label{induction 1}
\sum_{\ld\in\Par(d)}\log_p\vd_{p^r}(\ld)&=&\sum_{\ld\in\Par(d)}\log_p\vd_{p^{r-1}}(\ld)+
    \sum_{\ld\in\Par(d)}\left(\sum_{i=0}^{r-1}l(\ld^{(i)})+d_p(\ld^{(r-1)})\right),
\end{eqnarray}
where we interpret $\vd_1(\ld)=1$.
The Proposition follows easily by induction once we have proved the following equation
for all $r\geq 1$:
\[
\sum_{\ld\in\Par(d)}\left(\sum_{i=0}^{r-1}l(\ld^{(i)})+d_p(\ld^{(r-1)})\right)=l(d)\:.
\]
First note that each $\ld\in\Par(d)$
can uniquely be written as $\ld=\mu+p^r\mu'$, where
$\mu'\in\Par(j)$ and $\mu\in\Par_{p^r}(d-p^rj)$.
In fact, if
$\ld=\ld^{(0)}+p\ld^{(1)}+\cdots p^{N}\ld^{(N)}$ as before,
then $\mu=\sum_{i=0}^{r-1}p^i\ld^{(i)}\in\Par_{p^r}(d-p^rj)$;
note that for $\mu$, in the corresponding $p$-adic decomposition
we have $\mu^{(i)}=\ld^{(i)}$, for $i=0,\ldots,r-1$.
Hence,
\begin{eqnarray*}
\nonumber\sum_{\ld\in\Par(d)}\left(\sum_{i=0}^{r-1}l(\ld^{(i)})\right)&=&
\sum_{j\geq 0}p(j)\sum_{\mu\in\Par_{p^r}(d-p^rj)}l(\mu)\\
    &=&\sum_{j\geq 0}p(j)l_{p^r}(d-p^rj).
\end{eqnarray*}
Therefore, $\sum_{\ld\in\Par(d)}\sum_{i=0}^{r-1}l(\ld^{(i)})$ is the coefficient of $q^d$ in
\begin{eqnarray}\label{LHS 1}
P(q^{p^r})L_{p^r}(q)=P(q^{p^r}) P_{p^r}(q) T_{p^r}(q) = P(q)  T_{p^r}(q)\,,
\end{eqnarray}
where we have used formulae (\ref{formula:l-LPT}) and (\ref{formula:P_l}).

Obviously,
\begin{eqnarray*}
\sum_{\ld\in\Par(d)}d_p(\ld^{(r-1)})&=&
\sum_{j\geq0}p(j)\sum_{\mu\in\Par_{p^r}(d-p^rj)}d_p(\mu^{(r-1)})\:.
\end{eqnarray*}
Let $\tilde{c}_{p^r}(n)=\sum_{\mu\in\Par_{p^r}(n)}d_p(\mu^{(r-1)})$,
and $\widetilde{C}_{p^r}(q)=\sum_{n\geq0}\tilde{c}_{p^r}(n)q^n$.
Thus,  $\sum_{\ld\in\Par(d)}d_p(\ld^{(r-1)})$ is the coefficient of
$q^d$ in the generating series
\begin{eqnarray}\label{LHS 2}
P(q^{p^r})\widetilde{C}_{p^r}(q),
\end{eqnarray}
We now calculate the generating series $\widetilde{C}_{p^r}(q)$.
Note that for $\mu\in\Par_{p^r}(n)$, we may write $\mu=\eta+p^{r-1}\nu$,
where $\eta\in\Par_{p^{r-1}}(n-p^{r-1}j)$ and $\nu=\mu^{(r-1)}\in\Par_p(j)$.
Therefore,
\begin{eqnarray*}
\tilde{c}_{p^r}(n)&=&\sum_{j\geq0}p_{p^{r-1}}(n-p^{r-1}j)\sum_{\nu\in\Par_p(j)}d_p(\nu)\\
    &=&\sum_{j\geq0}p_{p^{r-1}}(n-p^{r-1}j)c_p(j),
\end{eqnarray*}
where for the second equation we have used (\ref{eqn: Cartan-defect}).
Hence, using Corollary~\ref{cor:Cartan-det} we obtain
\[
\widetilde{C}_{p^r}(q)=P_{p^{r-1}}(q)\Cg_p(q^{p^{r-1}})=\Cg_{p^r}(q).
\]
Now, by (\ref{LHS 2}) and Theorem~\ref{thm:Cartan-det}
we obtain that $\sum_{\ld\in\Par(d)}d_p(\ld^{(r-1)})$ is the coefficient of $q^d$ in
\begin{eqnarray}\label{LHS 3}
P(q^{p^r})\Cg_{p^r}(q)=P(q^{p^r})P_{p^r}(q)T(q^{p^r})=P(q)T(q^{p^r})\:.
\end{eqnarray}
Hence, adding (\ref{LHS 1}) and (\ref{LHS 3}), and using
formulae (\ref{formula:T}) and (\ref{formula:LPT}), we deduce that
\[
\sum_{\ld\in\Par(d)}\left(\sum_{i=0}^{r-1}l(\ld^{(i)})+d_p(\ld^{(r-1)})\right)
\]
is the coefficient of $q^d$ in
\[
P(q)T_{p^r}(q)+P(q)T(q^{p^r})=P(q)T(q)=L(q)\:.
\]
This proves the claim.
\end{pff}

\section{$\ell$-Blocks of Symmetric Groups}\label{kor}

In \cite{kor}, K\"{u}lshammer, Olsson, and Robinson developed the
theory of $\ell$-blocks of symmetric groups. The associated $\ell$-Cartan
matrix for $S_n$ is not unique. It depends on a choice of $\Z$-basis
for the $\Z$-span of the restriction of generalized characters of
$S_n$ to $\ell$-regular classes. Fix such a choice and define the
decomposition matrix $D_{\ell}(n)$ to be the transition matrix expressing
the restrictions of irreducible characters of $S_n$ to the
$\ell$-regular classes in terms of the characters in the fixed
$\Z$-basis. The $\ell$-Cartan matrix is then $C_{\ell}(n)=D_{\ell}(n)^tD_{\ell}(n)$.
Then it is shown in \cite{kor}, that two irreducible characters
belong to the same $\ell$-block exactly when the partitions labeling
them have the same $\ell$-core. The
determinant and invariant factors of the $\ell$-Cartan matrix of an
$\ell$-block depend only on the $\ell$-weight of the block, and not
on its $\ell$-core.  Donkin  has shown in  \cite{d}
that the invariant factors of this $\ell$-Cartan
matrix agree with those for the $\ell$-Cartan matrix of a block of the
Iwahori-Hecke algebra. In particular, it follows that the graded
invariant factors of the Cartan matrix for an $\ell$-block of the
symmetric group are given by the $\vd_{\ell}(\ld)$'s. For a positive
integer $k$, let again $\ell_k=\ell/(\ell,k)$,
and let $\pi_k$ be the set of primes
dividing $\ell_k$. For a partition $\mu\in\Par_{\ell}(n)$, define
\[
r_{\ell}(\mu)=
\prod_{k\geq 1}
\ell_k^{\lfloor \frac{m_k(\mu)}{\ell}\rfloor}\cdot\lfloor \frac{m_k(\mu)}{\ell}\rfloor!_{\pi_k}.
\]
We have the following ``KOR Conjecture'' (see \cite[Conjecture 6.4]{kor}):

\begin{cnj} The Cartan matrix $C_{\ell}(n)$ is unimodularly equivalent to a
diagonal matrix with entries $r_{\ell}(\mu)$ where $\mu$ runs through the set of
$\ell$-class regular partitions of $n$.
\end{cnj}

The goal of this section is to prove the following theorem:

\begin{thm}\label{kor conjecture new}
We have a  multiset equality
\[
\{r_{\ell}(\mu)\mid \mu\in\Par_{\ell}(n)\}
=
\{\vd_\ell(\ld)^{m_\ell^n(d)} \mid \ld\in \Par(d), d\leq \left\lfloor \frac n \ell \right\rfloor \}
\]
(the exponent in the second multiset is to be read as a multiplicity).
\end{thm}

This means that the explicit combinatorial descriptions for
a unimodularly equivalent diagonal form of the Cartan matrix
$C_\ell(n)$  given by the KOR Conjecture
and coming from Conjecture~\ref{conj}, respectively, coincide.

By Theorem~\ref{thm:main} we have thus the following contribution
towards the KOR conjecture:
\begin{cor}\label{kor conjecture}
Assume that in the prime decomposition
$\ell=\prod_{i=1}^rp_i^{r_i}$ we have $r_i\leq p_i$ for all
$i=1,\ldots,r$.
Then the KOR conjecture is true at $\ell$.
\end{cor}

Note also, that while the KOR conjecture is a conjecture about the full $\ell$-Cartan matrix
of~$S_n$, we may now also formulate a block version:

\begin{cnj}\label{conj:block Cartan}
Let $C_\ell(B)$ be the Cartan matrix of an $\ell$-block $B$ of $S_n$ of weight~$w$.
Then $C_{\ell}(B)$ is unimodularly equivalent to a
diagonal matrix with entries
\[
\vd_\ell(\ld)^{k(\ell-2, w-d)}\;, \ld\in \Par(d), d\leq w
\]
where again exponents are to be read as multiplicities.
\end{cnj}

Note that the generating function for $\sum_{d=0}^w p(d)k(\ell-2, w-d)$ is just
$P(q)^{\ell -1}$, and thus the size of the diagonal matrix is correct.
As evidence for the conjecture, we first confirm that also the determinant is correct.
Indeed, by the result on the determinant shown in the previous section,
we know that the product of all the numbers above (taking multiplicities into account)
is
\[
\ell^{\sum_{d=0}^w l(d)k(\ell -2,w-d)}\:.
\]
Now, by formula~(\ref{block-det}), we have $\sum_{d=0}^w l(d)k(\ell -2,w-d)=b_\ell(w)$,
and thus the conjectured diagonal matrix has indeed the correct determinant.

Furthermore, the largest number in the set  is
\[
\vd_\ell(1^w) = \ell^w w!_{\pi(\ell)}
\]
where $\pi(\ell)$ is the set of primes dividing $\ell$, and
$a_{\pi(\ell)}=\prod_{p\in \pi(\ell)} a_p$, for any integer~$a$.
In fact, it was shown in \cite[Theorems 6.1 and 6.2]{kor} that
$\ell^w w!_{\pi(\ell)}$ is the largest elementary divisor of $C_\ell(B)$.

\medskip

Towards the proof of Theorem~\ref{kor conjecture new}, we first want to collect the partitions
occurring for the two multisets into subsets associated with
a fixed $\ell$-class regular partition $\af$ of some $a\leq n$,
and then we will show that the corresponding contributed invariants
coincide for the two subsets.
\medskip

Let $\mu\in\Par_{\ell}(n)$.
As $\mu$ is $\ell$-class regular,
we can write it uniquely in the form $\mu=\hat{\mu}+\check{\mu}^\ell$,
where $\hat{\mu}$ is both $\ell$-class regular and $\ell$-regular,
and $\check{\mu}$ is $\ell$-class regular; here
$\check{\mu}^{\ell}=(1^{\ell m_1(\check{\mu})}2^{\ell
m_2(\check{\mu})}\ldots)$. Then clearly,
\[
r_{\ell}(\mu)=r_{\ell}(\check{\mu}^\ell)=\prod_{k\geq 1}\ell_k^{m_k(\check{\mu})}\cdot
(m_k(\check{\mu})!)_{\pi_k}.
\]
Also, given $\ld\in\Par(d)$, decompose $\ld=\ld^{(0)}+\ell\ld^{(1)}$,
where $\ld^{(0)}$ is $\ell$-class regular. Then, by definition,
\[
\vd_{\ell}(\ld)=
\vd_{\ell}(\ld^{(0)})=
\prod_{1\leq i\leq r}\prod_{\substack{k\geq 1\\0\leq\nu_p(k)<r_i}}
(p_i^{r_i})_k^{m_k(\ld^{(0)})}(m_k(\ld^{(0)})!)_{p_i}
\]
where $\ell=\prod_{i=1}^rp_i^{r_i}$ is the prime decomposition of $\ell$.

Now, let $m^n_{\ell}(d)$ be the
multiplicity of a graded invariant factor of degree $d$ in
$C_{\ell}(n)$ as in Theorem~\ref{total multiplicity},
i.e., any partition $\ld\in\Par(d)$ contributes $m^n_{\ell}(d)$ graded
invariant factors $\vd_{\ell}(\ld)$ to $C\ell(n)$.
We have the following lemma:

\begin{lem}
Let $a\leq n$ and $\af\in\Par_{\ell}(a)$. Then
$$
|\{\mu\in\Par_{\ell}(n)|\check{\mu}=\af\}|
=\sum_{d\geq1}m_{\ell}^n(d)|\{\ld\in\Par(d)|\ld^{(0)}=\af\}|.
$$
\end{lem}

\begin{pff}
Let LHS be the left hand side of the equation above, and RHS the right hand side.
We prove their equality by comparing their associated generating series.

First, observe that the LHS is the number of partitions of $n-\ell a$ that
are both $\ell$-regular and $\ell$-class regular. This is the coefficient
of $q^{n-\ell a}$ in the generating series (see section~2.2)
\[
\frac{P_{\ell}(q)}{P_{\ell}(q^\ell)}=\frac{P_{\ell}(q)}{P(q^\ell)}P(q^{\ell ^2}).
\]

We now turn to the right hand side. Observe that the $\ld$ appearing
there are partitions of $a+\ell j$, where $j\geq 0$. Counting
each such $\ld$ with its multiplicity $m_\ell^d(a+\ell j)$, we deduce that
\[
\mathrm{RHS}=\sum_{j\geq0}p(j)\sum_{w\geq
a+\ell j}d_\ell^0(n-\ell w)k(\ell-2,w-(a+\ell j)).
\]
Now,
\[
\sum_{w\geq a+\ell j}d^0_{\ell}(n-\ell w)k(\ell-2,w-(a+\ell j))
\]
is the coefficient of $q^{n-\ell (a+\ell j)}$ in the generating series
\[
\frac{P(q)}{P(q^\ell)^\ell}P(q^\ell)^{\ell -2}=\frac{P_{\ell}(q)}{P(q^\ell)}
\]
Hence, RHS is the coefficient of $q^{n-\ell a}$ in the generating series
\[
P(q^{\ell ^2})\frac{P_{\ell}(q)}{P(q^\ell)}
\]
which proves the lemma.
\end{pff}
\medskip

We know that for all partitions
$\mu\in\Par_{\ell}(n)$ with $\check{\mu}=\af$
we get the contribution
$r_\ell(\mu)=r_{\ell}(\af^\ell)$;
on the other hand, all partitions
$\ld\in\Par(d)$ with $\ld^{(0)}=\af$
give $m_\ell^n(d)$ contributions
$\vd_\ell(\ld)=\vd_{\ell}(\af)$.
Thus, Theorem~\ref{kor conjecture new} follows
from the next lemma:

\begin{lem}For $\af\in\Par_{\ell}$,
$$
r_{\ell}(\af^\ell)=\vd_{\ell}(\af).
$$
\end{lem}

\begin{pff} By definition, we have
\[
\vd_{\ell}(\af)=\prod_{1\leq i\leq r}\prod_{\substack{k\geq 1\\0\leq\nu_{p_i}(k)<r_i}}
(p_i^{r_i})_k^{m_k(\af)}(m_k(\af)!)_{p_i}.
\]
and
\[
r_{\ell}(\af^\ell)=\prod_{k\geq
1}\ell_k^{m_k(\af)}\cdot(m_k(\af)!)_{\pi_k}.
\]
Let $k$ be a part of $\af$. Since $\af\in\Par_{\ell}$, it follows
that $\ell\nmid k$. Write $k=\left(\prod_{i=1}^rp_i^{k_i}\right)k'$,
where $(\ell,k')=1$; note that there is at least one $j$ such that $k_j<r_j$.
Then,
\begin{eqnarray}\label{ell_k}
\ell_k=\prod_{i}p_i^{r_i-\min(r_ik_i)}
=\prod_{\substack{1\leq i\leq r\\k_i<r_i}}p_i^{r_i-k_i}
\end{eqnarray}
and, therefore,
\[
\prod_{k\geq 1}\ell_k^{m_k(\af)}=
\prod_{i=1}^r\prod_{\substack{k\geq1\\0\leq\nu_{p_i}(k)<r_i}}(p_i^{r_i})_k^{m_k(\af)}.
\]
Next, using (\ref{ell_k}), we deduce that $\pi_k=\{p_i|k_i<r_i\}$. Therefore,
\[
\prod_{k\geq 1}(m_k(\af)!)_{\pi_k}=\prod_{i=1}^r\prod_{0\leq\nu_{p_i}(k)<r_i}(m_k(\af)!)_{p_i}.
\]
\end{pff}

\bigskip

\noindent\textbf{Acknowledgements.}
We would like to thank the Mathematical Sciences Research Institute at Berkeley
for its hospitality and support
in the frame of the programs \emph{Combinatorial Representation Theory}
and \emph{Representation Theory of Finite Groups and Related Topics}.
We would also like to thank the program organizers for the invitation to take part
in these activities and for putting together these wonderful programs.

\end{document}